\documentclass[fleqn]{mat01}
\usepackage{times,mathtimy,amssymb,latexsym}
\begin{document}

\setcounter{page}{59}
\firstpage{59}

 \newcommand{\mbf}{\mathbb}
 \newcommand{\cl}{\mathcal}
 \newcommand{\ra}{\rightarrow}
 \newcommand{\BOX}{\rule{2mm}{2mm}}
 \newcommand{\ov}{\overline}
 \newcommand{\h}{\hat}
 \newcommand{\pa}{\partial}
 \newcommand{\sm}{\setminus}
 \newcommand{\la}{\lambda}
 \newcommand{\ti}{\times}
 \newcommand{\Si}{\Sigma}
 \newcommand{\ze}{\zeta}
 \newcommand{\Om}{\Omega}
 \newcommand{\til}{\tilde}
 \newcommand{\si}{\sigma}
 \newcommand{\ep}{\epsilon}
 \newcommand{\G}{\Gamma}

 \def\e{\mbox{\rm e}}

\newtheorem{theore}{Theorem}
\renewcommand\thetheore{\arabic{section}.\arabic{theore}}
\newtheorem{theor}[theore]{\bf Theorem}
\newtheorem{coro}[theore]{\rm COROLLARY}
\newtheorem{propo}[theore]{\rm PROPOSITION}
\newtheorem{lem}[theore]{Lemma}
\newtheorem{pot}[theore]{Proof of Theorem}
\newtheorem{case}{Case}
\def\remar{\trivlist \item[\hskip \labelsep{\it Remarks.}]}

\renewcommand\theequation{\thesection\arabic{equation}}


\title{Boundary regularity of correspondences in $\pmb{{\mbf C}^n}$}

\markboth{Rasul Shafikov and Kaushal Verma}{Boundary regularity of
correspondences in ${\mbf C}^n$}

\author{RASUL SHAFIKOV$^{1}$ and KAUSHAL VERMA$^{2}$}

\address{$^{1}$Department of Mathematics, Middlesex College, University
of Western Ontario, London, Ontario N6A 5B7, USA\\
\noindent $^{2}$Department of Mathematics, Indian Institute of
Science, Bangalore~560~012, India}

\volume{116}

\mon{February}

\parts{1}

\pubyear{2006}

\Date{MS received 6 October 2005}

\begin{abstract}
Let $M, M'$ be smooth, real analytic hypersurfaces of finite type
in ${\mbf C^n}$ and $\h f$ a holomorphic correspondence (not
necessarily proper) that is defined on one side of $M$, extends
continuously up to $M$ and maps $M$ to $M'$. It is shown that $\h
f$ must extend across $M$ as a locally proper holomorphic
correspondence. This is a version for correspondences of the
Diederich--Pinchuk extension result for CR maps.
\end{abstract}


\maketitle

\section{Introduction and statement of results}

\subsection{\it Boundary regularity}

Let $U, U'$ be domains in ${\mbf C}^n$ and let $M \subset U, M'
\subset U'$ be relatively closed, connected, smooth, real analytic
hypersurfaces of finite type (in the sense of D'Angelo). A~recent
result of Diederich and Pinchuk \cite{[DP3]} shows that a
continuous CR mapping $f\hbox{:}~M \ra M'$ is holomorphic in a
neighbourhood of $M$. The purpose of this note is to show that
their methods can be adapted to prove the following version of
their result for correspondences. We assume additionally that $M$
(resp. $M'$) divides the domain $U$ (resp. $U'$) into two
connected components $U^+$ and $U^-$ (resp. $U^{\prime\pm}$).

\begin{theor}[\!] Let ${\h f}\hbox{\rm :}~U^{-} \ra U'$ be a holomorphic
correspondence that extends continuously up to $M$ and maps $M$ to
$M'${\rm ,} i.e.{\rm ,} ${\h f}(M) \subset M'$. Then ${\h f}$
extends as a locally proper holomorphic correspondence across $M$.
\end{theor}
We recall that if ${\cl D} \subset {\mbf C}^p$ and ${\cl D'}
\subset {\mbf C}^m$ are bounded domains, a holomorphic
correspondence $\h f\hbox{:}~{\cl D} \ra {\cl D'}$ is a complex
analytic set $A \subset {\cl D} \times {\cl D'}$ of pure dimension
$p$ such that ${\ov A} \cap ({\cl D} \times {\pa {\cl D'}}) =
\emptyset$, where $\partial D'$ is the boundary of $D'$. In this
situation, the natural projection $\pi\hbox{:}~A \ra {\cl D}$ is
proper, surjective and a finite-to-one branched covering. If in
addition the other projection $\pi'\hbox{:}~A \ra {\cl D'}$ is
proper, the correspondence is called proper. The analytic set $A$
can be regarded as the graph of the multiple valued mapping ${\h
f} := \pi' \circ {\pi}^{-1}\hbox{:}~{\cl D} \ra {\cl D'}$. We also
use the notation $A=\text{Graph}(\hat f)$.

The branching locus $\sigma$ of the projection $\pi$ is a
codimension one analytic set in ${\cl D}$. Near each point in
${\cl D}\sm\sigma$, there are finitely many well-defined
holomorphic inverses of $\pi$. The symmetric functions of
these inverses are globally well-defined holomorphic functions on
${\cl D}$. To say that $\h f$ is continuous up to $\pa {\cl D}$
simply means that the symmetric functions extend continuously up
to $\pa {\cl D}$. Thus in Theorem~1.1 the various branches of $\h
f$ are continuous up to $M$ and each branch maps points on $M$ to
those on $M'$.

We say that $\h f$ in Theorem~1.1 extends as a {\it holomorphic
correspondence} across $M$ if there exist open neighbourhoods
$\tilde U$ of $M$ and $\tilde U'$ of $M'$, and an analytic set
$\tilde A\subset \tilde U\times \tilde U'$ of pure dimension $n$
such that (i)~${\rm Graph}(\h f)$ intersected with $(\tilde U \cap
U^-) \times (\tilde U'\cap U')$ is contained in $\tilde A$ and
(ii)~the projection $\tilde\pi\hbox{:}~\tilde A \to \tilde U$ is
proper. Without condition (ii), $\h f$ is said to extend as an
{\it analytic set}. Finally, the extension of $\hat f$ is a proper
holomorphic correspondence if in addition to (i) and (ii),
$\tilde\pi'\hbox{:}~\tilde A \to \tilde U'$ is also proper.

\setcounter{theore}{0}
\begin{coro}$\left.\right.$\vspace{.5pc}

\noindent Let $D$ and $D'$ be bounded pseudoconvex domains in
$\mathbb C^n$ with smooth real-analytic boundary. Let $\hat
f\hbox{\rm :}~D \to D'$ be a holomorphic correspondence. Then
$\hat f$ extends as a locally proper holomorphic correspondence to
a neighbourhood of the closure of $D$.
\end{coro}
The corollary follows immediately from Theorem~1.1 and \cite{[bs]}
where the continuity of $\hat f$ is proved. This generalizes a
well-known result of \cite{[br]} and \cite{[DF]} where the
extension past the boundary of $D$ is proved for holomorphic
mappings.

\subsection{\it Preservation of strata}

Let $M^+_s$ (resp. $M^-_s$) be the set of strongly pseudoconvex
(resp. pseudoconcave) points on $M$. The set of points where the
Levi form ${\cl L}_{\rho}$ has eigenvalues of both signs on
$T^{\mbf C}(M)$ and no zero eigenvalue will be denoted by
$M^{\pm}$ and finally $M^0$ will denote those points where ${\cl
L}_{\rho}$ has at least one zero eigenvalue on $T^{\mbf C}(M)$.
$M^0$ is a closed real analytic subset of $M$ of real dimension at
most $2n - 2$. Then
\begin{equation*}
M = M^+_s \cup M^-_s \cup M^{\pm} \cup M^0.
\end{equation*}
Further, let $M^+$ (resp. $M^-$) be the pseudoconvex (resp.
pseudoconcave) part of $M$, which equals the relative interior of
$\ov {M^+_s}$ (resp. $\ov {M^-_s}$). For non-negative integers $i,
j$ such that $i + j = n - 1$, let $M_{i, j}$ denote those points
at which ${\cl L}_{\rho}$ has exactly $i$ positive and $j$
negative eigenvalues on $T^{\mbf C}(M)$. Each (non-empty) $M_{i,
j}$ is relatively open in $M$ and semi-analytic whose relative
boundary is contained in $M^0$. With this notation, $M_{0, n - 1}
= M^-_s$ and $M_{n - 1, 0} = M^+_s$. Moreover, $M^{\pm}$ is the
union of all (non-empty) $M_{i, j}$ where both $i, j$ are at least
$1$ and $i + j = n - 1$. Note that points in $M^-_s, M^{\pm}$ are
in the envelope of holomorphy of $U^-$. Following \cite{[B]},
there is a semi-analytic stratification for $M^0$ given by
\begin{equation}
M^0 = \G_1 \cup \G_2 \cup \G_3 \cup \G_4,
\end{equation}
where $\G_4$ is a closed, real analytic set of dimension at most
$2n - 4$ and $\G_2 \cup \G_3 \cup \G_4$ is also a closed, real
analytic set of dimension at most $2n - 3$. Further, $\G_1, \G_2,
\G_3$ are either empty or smooth, real analytic manifolds; $\G_2,
\G_3$ have dimension $2n - 3$, and $\G_1$ has dimension $2n - 2$.
Finally, $\G_2$ and $\G_3$ are CR manifolds of complex dimension
$n - 2$ and $n - 3$ respectively. The set of points, denoted by
$\G^1_h$ in $\G_1$ where the complex tangent space to $\G_1$ has
dimension $n - 1$ is semi-analytic and has real dimension at most
$2n - 3$, as otherwise there would exist a germ of a complex
manifold in $M$ contradicting the finite type hypothesis. Then
\hbox{$\G_1\!\sm\!\G^1_h$} is a real analytic manifold of
dimension $2n - 2$ and has CR dimension $n - 2$. Using the same
letters to denote the various strata of $M^0$, there exists a
refinement of (1.1), so that $\G_1, \G_2, \G_3$ are all smooth,
real analytic manifolds of dimensions $2n - 2, 2n - 3, 2n - 3$
respectively, while the corresponding CR dimensions are $n - 2, n
- 2$, and $n - 3$. Finally, $\G_4$ is a closed, real analytic set
of dimension at most $2n - 4$.

\begin{theor}[\!]
With the hypothesis of Theorem~{\rm 1.1,} the extended
correspondence $\h f\hbox{\rm :}~M \ra M'$ satisfies the
additional properties{\rm :} $\h f(M^+) \subset {M}^{\prime
+}${\rm ,} $\h f(M^+ \cap M^0) \subset {M}^{\prime +} \cap
M^{\prime 0}$ and $\h f(M^-) \subset M^{\prime -}${\rm ,} $\h
f(M^- \cap M^0) \subset M^{\prime -} \cap M^{\prime 0}$.
Moreover{\rm ,} $\h f(M^+ \cap \G_j) \subset M^{\prime +} \cap
\G'_j$ and $\h f(M^- \cap \G_j) \subset M^{\prime -} \cap \G'_j$
for $j = 1, 3, 4$. Finally{\rm ,} $\h f$ maps the relative
interior of $\ov M^{\pm}$ to the relative interior of $\ov
M^{\prime\pm}$.
\end{theor}
Preservation of $\G_2$ is not always possible even for holomorphic
mappings as the following example shows: the domain $\Omega =
\{(z_1, z_2)\hbox{:}~\vert z_1 \vert^2 + \vert z_2 \vert^4 < 1\}$
is mapped to the unit ball in ${\mbf C}^2$ by the proper
holomorphic mapping $f(z_1, z_2) = (z_1, z^2_2)$. Points of the
form $\{(\e^{i \theta}, 0)\} \subset \partial \Omega$ are weakly
pseudoconvex and in fact form $\G_2 \subset \partial \Omega$, and
$f$ maps them to strongly pseudoconvex points.

\section{Segre varieties}

We will write $z = ('z, z_n) \in {\mbf C}^{n - 1} \times {\mbf C}$
for a point $z \in {\mbf C}^n$. The word `analytic' will always
mean complex analytic unless stated otherwise. The techniques of
Segre varieties will be used and here is a synopsis of the main
properties that will be needed. The proofs of these can be found
in \cite{[DF]} and \cite{[DW]}. As described above, let $M$ be a
smooth, real analytic hypersurface of finite type in ${\mbf C}^n$
that contains the origin. If $U$ is small enough, the
complexification $\rho(z, {\ov w})$ of $\rho$ is well-defined by
means of a convergent power series in $U \times U$. Note that
$\rho(z, {\ov w})$ is holomorphic in $z$ and anti-holomorphic in
$w$. For any $w \in U$, the associated Segre variety is defined as
\begin{equation*}
Q_w = \{z \in U\hbox{:}\ \rho(z, {\ov w}) = 0\}.
\end{equation*}
By the implicit function theorem, it is possible to choose
neighbourhoods $U_1 \subset \subset U_2$ of the origin such that
for any $w \in U_1$, $Q_w$ is a closed, complex hypersurface in
$U_2$ and
\begin{equation*}
Q_w = \{z = ('z, z_n) \in U_2\hbox{:}\ z_n = h('z, {\ov w}) \},
\end{equation*}
where $h('z, {\ov w})$ is holomorphic in $'z$ and anti-holomorphic
in $w$. Such neighbourhoods will be called a standard pair of
neighbourhoods and they can be chosen to be polydiscs centered at
the origin. It can be shown that $Q_w$ is independent of the
choice of $\rho$. For $\zeta \in Q_w$, the germ $Q_w$ at $\zeta$
will be denoted by ${}_{\zeta}Q_w$. Let ${\cl S} :=
\{Q_w\hbox{:}~w \in U_1\}$ be the set of all Segre varieties, and
let $\la\hbox{:}~w \mapsto Q_w$ be the so-called Segre map. Then
${\cl S}$ admits the structure of a finite dimensional analytic
set. It can be shown that the analytic\break set
\begin{equation*}
I_w := {\la}^{-1}(\la(w)) = \{z\hbox{:}\ Q_z = Q_w\}
\end{equation*}
is contained in $M$ if $w \in M$. Consequently, the finite type
assumption on $M$ forces $I_w$ to be a discrete set of points.
Thus $\la$ is proper in a small neighbourhood of each point of
$M$. For $w \in U_1^+$, the symmetric point ${}^sw$ is defined to
be the unique point of intersection of the complex normal to $M$
through $w$ and $Q_w$. The component of $Q_w \cap U_2^-$ that
contains the symmetric point is denoted by $Q_w^c$.

Finally, for all objects and notions considered above, we simply
add a prime to define their corresponding analogs in the target
space.

\section{Localization and extension across an open dense subset of $\pmb{M}$}

In the proof of Theorem~1.1 in order to show extension of $\h f$
as a holomorphic correspondence, it is enough to consider the
problem in an arbitrarily small neighbourhood of any point $p\in
M$. The reason is the following. Firstly, since the projection
$\pi\hbox{:}~\text{Graph}(\hat f) \to U^-$ is proper, the closure
of $\text{Graph}(\hat f)$ has empty intersection with $U^-\times
\partial U'$. Therefore, by \cite{[C]} \S~20.1, to prove the
continuation of $\hat f$ across $M$ as an analytic set, it is
enough to do that in a neighbourhood of any point in $M$.
Secondly, once the extension of $\hat f$ as a holomorphic
correspondence in a neighbourhood of any point $p\in M$ is
established, then globally there exists a holomorphic
correspondence defined in a neighbourhood $\tilde U$ of $M$ which
extends $\hat f$. To see that simply observe that if $F\subset
\tilde U\times \tilde{U'}$ is an analytic set extending $\h f$,
then by choosing smaller $\tilde U$ we may ensure that the
projection to the first component is proper, as otherwise there
would exist a point $z$ on $M$ such that $\hat F(z)$ has positive
dimension (here $\hat F$ is the map associated with the set $F$).
This however contradicts local extension of $\hat f$ near $z$ as a
holomorphic correspondence.

Since the projection $\pi\hbox{:}~\text{Graph}(\hat f) \to U^-$ is
proper, ${\rm Graph}(\h f)$ is contained in the analytic set $A
\subset U^- \times U'$, defined by the zero locus of holomorphic
functions $P_1(z, z'_1), P_2(z, z'_2), \dots, P_n(z, z'_n)$ given
by \setcounter{equation}{0}
\begin{equation}\label{a}
P_j(z, z'_j) = {z'_j}^l + a_{j1}(z) {z'_j}^{l - 1} + \cdots +
a_{jl}(z),
\end{equation}
where $l$ is the generic number of images in $\hat f(z)$, and $1
\leq j \leq n$ (for details, see \cite{[C]}). The coefficients
$a_{\mu \nu}(z)$ are holomorphic in $U^-$ and extend continuously
up to $M$. This is the definition of continuity of the
correspondence $\hat f$ up to $M$ which is equivalent to that
given in \S1.1.

The discriminant locus is $\{R_j(z) = 0\}$, $1 \leq j \leq n$,
where $R_j(z)$ is a universal polynomial function of $a_{j
\mu}(z)$ ($1 \leq \mu \leq l$) and hence by the uniqueness
theorem, it follows that ${\ov {\{R_j(z) = 0\}}} \cap M$ is
nowhere dense in $M$, for all $j$. The set of points $S$ on $M$
which do not belong to ${\ov {\{R_j(z) = 0\}}} \cap M$ for any $j$
is therefore open and dense in $M$. Near each point $p$ on $S$,
$\h f$ splits into well-defined holomorphic maps $f_1(z), f_2(z),
\dots, f_l(z)$ each of which is continuous up to $M$.

If $p \in S \cap (M^- \cup M^{\pm})$, the functions $a_{\mu
\nu}(z)$ extend holomorphically to a neighbourhood of $p$ and
hence $\h f$ extends as a holomorphic correspondence across $p$.
It is therefore sufficient to show that $\h f$ extends across an
open dense subset of $S \cap M^+$. But this follows from Lemma~3.2
and Corollary~3.3 in \cite{[DP3]}. We denote by $\Si \subset M$
the non-empty open dense subset of $M$ across which $\h f$ extends
as a holomorphic correspondence.

\section{Extension as an analytic set}

Fix $0 \in M$ and let $p'_1, p'_2, \dots, p'_k \in \h f(0) \cap
M'$. The continuity of $\h f$ allows us to choose neighbourhoods
$0 \in U_1$ and $p'_i \in U'_i$ and local correspondences $\h
f_i\hbox{:}~U_1^- \ra U'_i$ that are irreducible and extend
continuously up to $M$. Moreover, $\h f_i(0) = p'_i$ for all $1
\leq i \leq k$. It will suffice to focus on one of the $\h
f'_i$'s, say $\h f'_1$ and to show that it extends holomorphically
across the origin. Abusing notation, we will write $\h f'_1 = \h
f$, $U'_1 = U'$ and $p'_1 = 0'$ . Thus $\h f\hbox{:}~U_1^- \ra U'$
is an irreducible holomorphic correspondence and $\h f(0) = 0'$.
Define
\begin{equation*}
V^+ = \{ (w, w') \in U_1^+ \ti U'\hbox{:}\ \h f(Q_w^c) \subset
Q'_{w'} \}.
\end{equation*}
Then $V^+$ is non-empty. Indeed, $\h f$ extends across an open
dense set near the origin and \cite{[V]} shows that the invariance
property of Segre varieties then holds. Moreover, a similar
argument as in \cite{[S2]} shows that $V^+ \subset U_1^+ \ti U'$
is an analytic set of dimension $n$ and exactly the same arguments
as in Lemmas~4.2~-- 4.4 of \cite{[DP3]} show that: first, the
projection $\pi\hbox{:}~V^+ \ra U^+ := \pi(V^+) \subset U_1^+$ is
proper (and hence that $U^+ \subset U_1^+$ is open) and second,
the projection $\pi'\hbox{:}~V^+ \ra U'$ is locally proper. Thus,
to $V^+$ is associated a correspondence $F^+\hbox{:}~U^+ \ra U'$
whose branches are $\h F^+ = \pi' \circ \pi^{-1}$.

Let $a \in M$ be a point close to the origin, across which $\h f$
extends as a holomorphic correspondence. If $\h f$ is well-defined
in the ball $B(a, r)$, $r > 0$ and $w \in B(a, r)^-$, it follows
from Theorem~4.1 in \cite{[V]} that all points in $\h f(w)$ have
the same Segre variety. By analytic continuation, the same holds
for all $w \in U_1^-$. Using this observation, it is possible to
define another correspondence $F^-\hbox{:}~U_1^- \ra U'$ whose
branches are $\h F^-(w) = (\la')^{-1} \circ {\la'} \circ \h f(w)$.
Let $U := U_1^- \cup U^+ \cup (\Si \cap U_1)$. The invariance
property of Segre varieties shows that the correspondences $\h
F^+, \h F^-$ can be glued together near points on $\Si \cap U_1$.
Hence, there is a well-defined correspondence $\h F\hbox{:}~U \ra
U'$ whose values over $U^+$ and $U_1^-$ are $\h F^+$ and $\h F^-$
respectively. Note that
\begin{equation*}
F := {\rm Graph}(\h F) = \{(w, w') \in U \ti U'\hbox{:}\ w' \in \h F(w)
\}
\end{equation*}
is an analytic set in $U \ti U'$ of pure dimension $n$, with
proper projection $\pi\hbox{:}~F \ra U$. Once again, the
invariance property shows that all points in $\h F(w)$, $w \in U$,
have the same Segre variety.

\setcounter{theore}{0}
\begin{lem}\label{l1}
The correspondence $\h F$ satisfies the following properties{\rm
:}
\begin{enumerate}
\renewcommand\labelenumi{\rm (\roman{enumi})}
\leftskip .3pc
\item For $w_0 \in \pa U \cap U_1^+${\rm ,} ${\rm cl}_{\h F}(w_0) \subset
\pa U'$.

\item ${\rm cl}_{\h F}(0) \subset Q'_{0'}$.

\item If ${\rm cl}_{\h F}(0) = \{0'\}${\rm ,} then $0 \in \Si$.

\item $F \subset (U_1\!\sm\!(M\!\sm\!\Si)) \ti U'$ is a closed analytic
set.\vspace{-.5pc}
\end{enumerate}
\end{lem}

\begin{proof}$\left.\right.$\vspace{-.5pc}

\begin{enumerate}
\renewcommand\labelenumi{\rm (\roman{enumi})}
\leftskip .3pc \item Choose $(w_j, w'_j) \in F$ converging to
$(w_0, w'_0) \in (\pa U \cap U_1^+) \ti {\ov U}'$. Then $\h
f(Q^c_{w_j}) \subset Q'_{w'_j}$ for all $j$. If $w'_0 \in U'$,
then passing to the limit, we get $\h f(Q^c_{w_0}) \subset
Q'_{w'_0}$ which shows that $(w_0, w'_0) \in F$ and hence $w_0 \in
U$, which is a contradiction. This also\break proves (iv).

\item Choose $w_j \in U$ converging to $0$. There are two cases to
consider. First, if $w_j \in U_1^- \cup (\Si \cap U_1)$ for all
$j$, it follows that $\h f(w_j) \ra 0'$. Moreover, for any $w'_j
\in \h F(w_j)$, $Q'_{w'_j} = Q'_{\h f(w_j)}$. If $U'$ is small
enough, the equality $Q'_{w'} = Q'_{0'}$ implies that $w' = 0'$
and thus we conclude that $w'_j \ra 0' \in Q'_{0'}$. Second, if
$w_j \in U^+$ for all $j$, then $\h f(Q^c_{w_j}) \subset
Q'_{w'_j}$ for any $w'_j \in \h F(w_j)$. Let $w'_j \ra w'_0 \in
U'$. If $\ze \in Q^c_{w_j}$, then $\h f(\ze) \in Q'_{w'_j} \ra
Q'_{w'_0}$. But $w_j \ra 0$ implies that ${\rm dist}(Q^c_{w_j}, 0)
\ra 0$ and hence $\h f(\ze) \ra 0'$. Thus $0' \in Q'_{w'_0}$ which
shows that $w'_0 \in Q'_{0'}$.

\item If ${\rm cl}_{\h F}(0) = \{0'\}$, then (i) shows that $0
\notin \pa U \cap U_1^+$. Let $B(0, r)$ be a small ball around the
origin such that $B(0, r) \cap \pa U = \emptyset$. The
correspondence $\h F$ over $B(0, r)^+$ is the union of some
components of the zero locus of a system of monic
pseudo-polynomials whose coefficients are bounded holomorphic
functions on $B(0, r)^+$. By Trepreau's theorem, all these
coefficients extend holomorphically to $B(0, r)$, and the extended
zero locus contains the graph of $\h f$ near the origin since
$\Si$ is dense. It follows that $0 \in \Si$. \hfill
$\Box$\vspace{-2.2pc}
\end{enumerate}
\end{proof}
Following \cite{[S1]}, for any $w_0 \in U$, it is possible to find
a neighbourhood $\Omega$ of $w_0$, relatively compact in $U$ and a
neighbourhood $V \subset U_1$ of $Q_{w_0} \cap U_1$ such that for
$z \in V$, $Q_z \cap \Omega$ is non-empty and connected.
Associated with the pair $(\Omega, V)$ is \setcounter{equation}{0}
\begin{equation}\label{tildef}
\til F := \til F(w_0, \Om, V) = \{(z, z') \in V \ti U' : \h F(Q_z
\cap \Om) \subset Q'_{z'} \}
\end{equation}
which (see \cite{[DP4]}) is an analytic set of dimension at most
$n$. If $w_0 \in \Si$, then Corollary~5.3 of \cite{[DP3]}, shows
that $F \cap (V \ti U')$ is the union of irreducible components of
$\til F$ of dimension $n$. As in \cite{[DP3]} we call $(w_0,z_0)
\in U\times Q_{w_0}$ a pair of reflection if there exist
neighbourhoods $\Omega(w_0)\ni w_0$ and $\Omega(z_0)\ni z_0$ such
that for all $w\in\Omega(w_0)$, $\hat F (Q_w\cap
\Omega(z^0))\subset Q'_{\hat F(w)}$. It follows from the
invariance property of Segre varieties that the definition of the
pair of reflection is symmetric. As an example we note that if the
set $\tilde F$ defined in \eqref{tildef} contains $F \cap (V \ti
U')$, then $(w_0,z)$ is a point of reflection for any point $z$ in
a connected component of $Q_{w_0}\cap U$ containing $w_0$.

Let $w_0 \in U$, $z_0 \in Q_{w_0} \cap \Si$ be a pair of
reflection. Fix $B(z_0, r)$, a small ball around $z_0$ where $\h
f$ is well-defined and let $S(w_0, z_0) \subset \til F \cap
((Q_{w_0} \cap U_1) \ti U')$ be the union of those irreducible
components that contain ${\rm Graph}(\h f)$ over $Q_{w_0} \cap
B(z_0, r)$. Note that $S(w_0, z_0)$ is an analytic set of
dimension $n - 1$ and is contained in $(Q_{w_0} \cap U_1) \ti U'$
and moreover, the invariance property shows that
\begin{equation*}
S(w_0, z_0) \subset ( (Q_{w_0} \cap U_1) \ti (Q'_{\h F(w_0)} \cap
U')).
\end{equation*}
Furthermore, from the above considerations it follows that for any
$z\in \pi(S(w_0,z_0))$ the point $(w_0,z)$ is a pair of
reflection. Finally, let the cluster set of a sequence of closed
sets $\{C_j\} \subset {\mathcal D}$, where ${\mathcal D}$ is some
domain, be the set of all possible accumulation points in
${\mathcal D}$ of all possible sequences $\{c_j\}$ where $c_j \in
C_j$.\vspace{-.2pc}

\setcounter{theore}{0}
\begin{propo}\label{p}$\left.\right.$\vspace{.5pc}

\noindent Let $\{z_{\nu}\} \in \Si$ converge to $0$. Suppose that
the cluster set of the sequence $\{S(z_{\nu}, z_{\nu})\}$ contains
a point $(\ze_0, \ze'_0) \in U \ti U'$. Then $\h f$ extends as an
analytic set across the origin.\vspace{-.2pc}
\end{propo}

\begin{proof} First, the pair $(z_{\nu}, z_{\nu})$ is an example of a
pair of reflection and hence $S(z_{\nu}, z_{\nu})$ is
well-defined. Also, note that $(z_{\nu}, \h f(z_{\nu})) \ra (0,
0')$. Choose $(\ze_{\nu}, \ze'_{\nu}) \in S(z_{\nu}, z_{\nu})$
that converges to $(\ze_0, \ze'_0) \in U \ti U'$. It follows that
$(\ze_{\nu}, z_{\nu})$ is a pair of reflection. Let $\Om, V$ be
neighbourhoods of $\ze_0$ and $Q_{\ze_0}$ as in the definition of
$\til F(\ze_0, \Om, V)$. Since $\ze_0 \in U$, it follows that
$\til F(\ze_0, \Om, V)$ is a non-empty, analytic set in $V \times
U'$. Shrinking $U_1$ if needed, $Q_{\ze_{\nu}} \cap U_1 \subset V$
and $\ze_{\nu} \in \Om$ for all large $\nu$. This shows that $\til
F(\ze_{\nu}, \Om, V) = \til F(\ze_0, \Om, V)$ for all large $\nu$.
Lemma~5.2 of \cite{[DP3]} shows that $\til F(\ze_{\nu}, \Om, V)$
contains the graph of all branches of $\h f$ near $z_{\nu}$ and
hence $\til F(\ze_0, \Om, V)$ contains the graph of $\h f$ near
$(0, 0')$. Therefore, $\til F(\ze_0, \Om, V)$ extends the graph of
$\h f$ across the origin. \hfill $\Box$\vspace{-.2pc}
\end{proof}

\begin{remar}
First, as in \cite{[DP3]} this proposition will be valid if the
pair $(z_{\nu}, z_{\nu})$ were replaced by a pair of reflection
$(w_{\nu}, z_{\nu}) \in U \ti \Si$ that converges to $(0, 0')$ and
$\hat F(w_\nu)$ \hbox{clusters} at some point in $U'$. Second, this
proposition shows the relevance of studying the \hbox{cluster} set of a
sequence of analytic sets (see \cite{[SV]} also). In general, the
hypothesis that the cluster set of $\{ S(z_{\nu}, z_{\nu}) \}$ (or
$S(w_{\nu}, z_{\nu})$ in case $(w_{\nu}, z_{\nu})$ is a pair of
reflection) contains a point in $U \ti U'$ cannot be guaranteed
since the projection $\pi\hbox{:}~S(z_{\nu}, z_{\nu}) \ra U$ is
not known to be proper. However, the following version of
Lemma~5.9 in \cite{[DP3]}\break holds.
\end{remar}

\begin{lem}\label{l2}
There are sequences $(w_{\nu}, z_{\nu}) \in U \ti \Si$, $w'_{\nu}
\in \h F(w_{\nu})$ and analytic sets $\si_{\nu} \subset U$ of pure
dimension $p \geq 1$ {\rm (}$p$ independent of $\nu${\rm )} such
that{\rm :}\vspace{-.3pc}

\begin{enumerate}
\renewcommand\labelenumi{\rm (\roman{enumi})}
\leftskip .2pc
\item $(w_{\nu}, z_{\nu}) \ra (0, 0)$ and $(w_{\nu}, z_{\nu})$ is a
pair of reflection for all $\nu$.

\item $w'_{\nu} \ra w'_0 \in U'$ and $z_{\nu} \in \si_{\nu} \subset
\pi(S(w_{\nu}, z_{\nu}))$.\vspace{-.6pc}
\end{enumerate}
\end{lem}

\begin{proof}
Choose a sequence $z_{\nu} \in \Si$ that converges to the origin.
If the projections $\pi\hbox{:}~S(z_{\nu}, z_{\nu}) \ra U$ were
proper for all $\nu$, then for some fixed $r > 0$ and $\nu$ large
enough, let $\si_{\nu} := Q_{z_{\nu}} \cap B(z_{\nu}, r)$,
$w_{\nu} = z_{\nu}$ and $w'_{\nu} \in \h f(z_{\nu})$. It can be
seen that the lemma holds with these choices. On the other hand,
if $\pi$ is not known to be proper on $S(z_{\nu}, z_{\nu})$, no
fixed value of $r$, as described above, exists. Hence, for
arbitrarily small values of $r' > 0$, there exist $(w_{\nu},
w'_{\nu}) \in S(z_{\nu}, z_{\nu}) \cap (U^+ \ti U')$ such that
$w_{\nu} \ra 0$ and $w'_{\nu} \ra w'_0$ with $\vert w'_0 \vert =
r'$. Since $M'$ is of finite type, we may assume that $Q'_{w'_0}
\not= Q'_{0'}$. Moreover, note that $w'_0 \in Q'_{0'} \cap U'$
(which shows that $0' \in Q'_{w'_0}$) and $(w_{\nu}, z_{\nu})$ is
a pair of reflection for all $\nu$. By making a small holomorphic
perturbation of coordinates in the target space, if needed, it
follows that $0' \in Q'_{w'_0} \cap \{z' \in U'\hbox{:}~z'_2 =
\cdots = z'_n = 0 \}$ is an isolated point. Therefore, there
exists an $\ep > 0$ such that after shrinking $U'$, if needed,
$q'_0 := Q'_{w'_0} \cap \{ z' \in U'\hbox{:}~z'_2 = \cdots = z'_{n
- 1} = 0, \vert z'_n \vert < \ep \}$ (which is an analytic set of
dimension $1$ in $U' \cap \{ \vert z'_n \vert < \ep \}$ containing
the origin) has no limit points on $\pa U' \cap \{ \vert z'_n
\vert < \ep \}$. Let $l$ be the multiplicity of $\h
f\hbox{:}~U_1^- \ra U'$. Let $\h f(z_{\nu}) = \{ \ze_{\nu}^j \}, 1
\leq j \leq l$ counted with multiplicity. For large $\nu$, the $l$
sets
\begin{equation*}
q'_{\nu, j} = Q'_{w'_{\nu}} \cap \{z' \in U'\hbox{:}\  z'_k =
{(\ze_{\nu}^j)}_k, \quad 2 \leq k \leq n - 1, \quad \vert z'_n
\vert < \ep \}
\end{equation*}
are analytic, of dimension $1$, in $U' \cap \{ \vert z'_n \vert <
\ep \}$ without limit points on $\pa U' \cap \{ \vert z'_n \vert <
\ep \}$ and clearly contain $(z_{\nu}, \ze_{\nu}^j)$. Since
$\pi'(S(w_{\nu}, z_{\nu})) \subset Q'_{w'_{\nu}}$,
\begin{equation*}
s_{\nu, j} := S(w_{\nu}, z_{\nu}) \cap \{(z, z')\hbox{:}\  z'_k =
{(\ze_{\nu}^j)}_k , \quad 2 \leq k \leq n - 1 \}
\end{equation*}
 are analytic
sets of dimension at least $1$ in $U_1 \ti (U' \cap \{ \vert z'_n
\vert < \ep \})$ for all $1 \leq j \leq l$. By construction, the
analytic sets $q'_{\nu, j}$ do not have limit points on $\pa U'
\cap \{ \vert z'_n \vert < \ep \}$ and hence $s_{\nu, j}$ do not
have limit points on $U_1 \ti (\pa U' \cap \{ \vert z'_n \vert <
\ep \})$. By Lemma~\ref{l1}, ${\rm cl}_{\h F}(0) \subset Q'_{0'} =
\{ z'_n = 0 \}$ and by shrinking $U_1$ if needed, this shows that
$s_{\nu, j}$ have no limit points on $U_1 \ti (U' \cap \{ \vert
z'_n \vert = \ep \})$. Thus for large $\nu$ and all $j$, the
projections $\pi\hbox{:}~s_{\nu, j} \ra U_1$ are proper and their
images $\si_{\nu, j} := \pi(s_{\nu, j})$ are analytic sets in
$U_1$ of dimension at least $1$ and $z_{\nu} \in \si_{\nu, j}$ for
all $\nu, j$. It remains to pass to subsequences if necessary to
choose $\si_{\nu, j}$ with constant dimension. \hfill $\Box$
\end{proof}

One conclusion that follows now is: if $\h f$ does not extend as
an analytic set across the origin, then ${\rm cl}(\si_{\nu})
\subset M\!\sm\!\Si$. Indeed, if there exists $\ze_0 \in {\rm
cl}(\si_{\nu}) \cap (U_1\!\sm\!(M\!\sm\!\Si))$, let $(\ze_{\nu},
\ze'_{\nu}) \in S(w_{\nu}, z_{\nu})$ converge to $(\ze_0, \ze'_0)
\in U_1 \ti U'$. Proposition~\ref{p} now shows that $\ze_0 \in \pa
U \cap U_1$. But since $\ze_0 \notin M\!\sm\!\Si$, it follows from
Lemma~\ref{l1} that $\ze'_0 \in \pa U'$ which is a contradiction.

The goal will now be to show that $\h f$ extends as an analytic
set across the origin. For this, choose $\{ z_{\nu} \} \in \Si$
converging to the origin and consider the analytic sets
$S(z_{\nu}, z_{\nu})$. By Proposition~\ref{p}, it suffices to show
that $\pi({\rm cl}(S(z_{\nu}, z_{\nu})) \cap U \not = \emptyset$.
Let
\begin{equation*}
S' := \pi' ( {\rm cl}(S(z_{\nu}, z_{\nu})) \cap (\{0\} \ti U'))
\subset Q'_{0'}
\end{equation*}
and let $m$ be the dimension of $\h S'$~--~the smallest closed
analytic set containing $S'$ (the so-called Segre completion of
\cite{[DP3]}). If $m = 0$, then $0'$ is an isolated point in $S'$
and after shrinking $U_1, U'$ suitably, it follows that ${\rm
cl}(S(z_{\nu}, z_{\nu}))$ has no limit points on $U_1 \ti \pa U'$.
Thus $\pi\hbox{:}~S(z_{\nu}, z_{\nu}) \ra U_1$ are proper
projections and therefore $\pi(S(z_{\nu}, z_{\nu})) = Q_{z_{\nu}}
\cap U_1$ for all large $\nu$. Hence $\pi({\rm cl}(S(z_{\nu},
z_{\nu}))) = Q_0 \cap U_1$. If $\h f$ did not extend as an
analytic set across the origin, the aforementioned remark shows
that with $\si := Q_{z_{\nu}} \cap U_1$, $Q_0 \cap U_1 = {\rm
cl}(\si_{\nu}) \subset M\!\sm\!\Si \subset M$. This cannot happen
as $M$ is of finite type. Hence $\h f$ extends as an analytic set
across the origin in case $m = 0$. We may therefore suppose that
$m > 0$. We recall the following lemma proved by Diederich and
Pinchuk:

\begin{lem}\hskip -.3pc{\rm (\cite{[DP3]},
{\it Lemma}~9.8).}\label{dpl}\ \ \ Let $S'$ be a subset of $Q'_{0'}${\rm
,} $0'\in S'$ and $m=\dim \hat S'$. Then after possibly shrinking
$U_1${\rm ,} there are points $w^{\prime 1},\dots,w^{\prime k}\in
S'$ $(k\leq n-1)$ such that one of the following holds{\rm :}

\begin{enumerate}
\renewcommand\labelenumi{\rm (\arabic{enumi})}
\leftskip .15pc
\item $k=m$ and $\dim(\hat S'\cap Q'_{w^{\prime 1}}\cap\cdots\cap
Q'_{w^{\prime k}})=0${\rm ;}

\item $k\geq 2m-n+1$ and $\dim(\hat S'\cap Q'_{w^{\prime
1}}\cap\cdots\cap Q'_{w^{\prime k}})=m-k$.\vspace{-.5pc}
\end{enumerate}
\end{lem}
Thus there are two cases to consider.

\begin{case}{\rm
Choose $(w_{1 \nu}, w'_{1 \nu}), (w_{2 \nu}, w'_{2 \nu}), \dots,
(w_{m \nu}, w'_{m \nu}) \in S(z_{\nu}, z_{\nu})$ so that $w_{\mu
\nu} \ra 0$ and $w'_{\mu \nu} \ra w'_{\mu}$ for all $1 \leq \mu
\leq m$. A~generic choice of $w_{\mu \nu}$ (see p.~136 in
\cite{[DP3]}) ensures that $q_{m \nu} := Q_{w_{1 \nu}} \cap
Q_{w_{2 \nu}} \cap \cdots \cap Q_{w_{m \nu}}$ has dimension $n -
m$. Each $(w_{\mu \nu}, z_{\nu})$ is a pair of reflection and
hence the analytic set
\begin{equation*}
S^m_{\nu} := \bigcap_{1 \leq \mu \leq m} S(w_{\mu \nu}, z_{\nu})
\subset (q^{m \nu} \ti {q}^{\prime m \nu}) \cap (U_1 \ti U')
\end{equation*}
is well-defined. If $m = n - 1$, then Lemma~9.7 of \cite{[DP3]}
shows that the germ of ${q}^{\prime(n - 1)}$ at the origin has
dimension 1. Moreover, $\h S' = Q'_{0'}$ and Lemma~\ref{dpl}
implies that ${q}^{\prime(n - 1)} \cap Q'_{0'}$ contains $0'$ as
an isolated point. Since ${\rm cl}_{\h F}(0) \subset Q'_{0'}$, it
follows that $0'$ is an isolated\break point of
\begin{equation*}
\pi'({\rm cl}(S^{n - 1}_{\nu}) \cap (\{0\} \ti U')) \subset
{q}^{\prime(n - 1)} \cap Q'_{0'} = \{0'\}.
\end{equation*}
Shrinking $U_1$, the projection $\pi\hbox{:}~S^{n - 1}_{\nu} \ra
U_1$ becomes proper and $\pi(S^{n - 1}_{\nu}) = q^{n - 1, \nu}
\cap U_1$. By Theorem~7.4 of \cite{[DP3]}, there is a subsequence
of $q^{n - 1, \nu} \cap U_1$ that converges to an analytic set $A
\subset U_1$ of pure dimension 1 and contains the origin. $A$
contains points $\ze_0$ that do not belong to $M$ because of the
finite type assumption and $\ze_0 \in \pi({\rm cl}(S^{n -
1}_{\nu})) \subset \pi({\rm cl}(S(w_{\mu \nu}, z_{\nu})))$. By
Proposition~\ref{p}, $\h f$ extends as an analytic set across the
origin.

If $m < n - 1$, the dimension of $S^m_{\nu} \cap S(z_{\nu},
z_{\nu})$ is at least $n - m - 1 > 0$. Now
\begin{equation*}
\pi'({\rm cl}(S^m_{\nu} \cap S(z_{\nu}, z_{\nu})) \cap (\{0\} \ti
U')) \subset {q}^{\prime m} \cap {\h S'} = \{0'\},
\end{equation*}
the last equality following from Lemma~\ref{dpl}. The projection
$\pi\hbox{:}~S^m_{\nu} \cap S(z_{\nu}, z_{\nu}) \ra U_1$ is
therefore proper for small $U_1$ and that $\pi(S^m_{\nu} \cap
S(z_{\nu}, z_{\nu})) = q^{m \nu} \cap Q_{z_{\nu}} \cap U_1$.
Again, by Theorem~7.4 of \cite{[DP3]}, there is a subsequence of
$q^{m \nu} \cap Q_{z_{\nu}} \cap U_1$ that converges to an
analytic set $A \subset U_1$ of positive dimension and as before
this shows that $\h f$ extends as an analytic set across the
origin.}
\end{case}

\begin{case}
{\rm As before, choose $(w_{1 \nu}, w'_{1 \nu}), (w_{2 \nu}, w'_{2
\nu}), \dots, (w_{k \nu}, w'_{k \nu}) \in S(z_{\nu}, z_{\nu})$
such that $w_{\mu \nu} \ra 0$ and $w'_{\mu \nu} \ra w'_{\mu}$ for
all $1 \leq \mu \leq k$ and $q_{k \nu} = Q_{w_{1 \nu}} \cap
Q_{w_{2 \nu}} \cap \cdots \cap Q_{w_{k \nu}}$, ${\tilde q}^{k \nu}
:= Q_{z_{\nu}} \cap q^{k \nu}$ have dimension $n - k$ and $n - k -
1$ respectively. Now note that ${\rm dim}(S^k_{\nu} \cap
S(z_{\nu}, z_{\nu})) \geq n - k - 1 > 1$. Indeed, the inequalities
$2m - n + 1 \leq k < m$ show that $m \leq n - 2$ and hence $k < n
- 2$. Since the dimension of $\h S' \cap {q}^{\prime k}$ is $m -
k$, choose coordinates so that
\begin{equation*}
\h S' \cap {q}^{\prime k} \cap \{z' \in U'\hbox{:}\  z'_1 = z'_2 = \cdots
= z'_{m - k} = 0\} = \{0'\}.
\end{equation*}
Let $\h f(z_{\nu}) = \{\ze^j_{\nu}\}$, $1 \leq j \leq l$, $l$
being the multiplicity of $\h f$. The $l$ sets
\begin{align*}
T_{\nu, j} &= \{(z, z') \in S^k_{\nu} \cap S(z_{\nu}, z_{\nu})\hbox{:}\
z'_1 = {(\ze^j_{\nu})}_1,\\[.3pc]
&\quad\ z'_2 = {(\ze^j_{\nu})}_2, \dots, z'_{m - k} =
{(\ze^j_{\nu})}_{m - k}\},
\end{align*}
where $1 \leq j \leq l$ are analytic sets in $U_1 \ti U'$ and have
dimension at least $n - k - 1 - (m - k) = n - m - 1 > 0$. By
construction,
\begin{align*}
&\pi' ({\rm cl}(T_{\nu, j}) \cap ( \{0\} \ti U')) \subset \h S'
\cap {q}^{\prime k}\\[.3pc]
&\quad\ \cap \{z' \in U'\hbox{:}\ z'_1 = z'_2 = \cdots = z'_{m - k} = 0\}
= \{0'\}
\end{align*}
and hence by shrinking $U_1, U'$, the projections
$\pi\hbox{:}~T_{\nu, j} \ra U_1$ are proper and the images
$\si_{\nu, j} := \pi(T_{\nu, j}) \subset U_1$ are analytic and
have dimension $n - m - 1$. Moreover $\si_{\nu, j} \subset {\tilde
q}^{k \nu}$, and since ${\tilde q}^{k \nu}$ depend
anti-holomorphically on the $k$-tuple defining it, Theorem~7.4 of
\cite{[DP3]} shows that ${\tilde q}^{k \nu}$ converges to an
analytic set $\tilde A \subset U_1$ of dimension $n - k - 1$,
after passing to a subsequence. Working with this subsequence, we
see that ${\rm cl}(\si_{\nu, j}) \subset \tilde A$. On the other
hand, since $2m - n + 1 \leq k$, it follows, as in \cite{[DP3]},
that
\begin{equation*}
{\rm dim}\; {\tilde A} = n - k - 1 \leq 2(n - m - 1) = 2 \; {\rm
dim}\;{\si_{\nu, j}}.
\end{equation*}
Proposition~8.3 of \cite{[DP3]} shows that ${\rm cl}(\si_{\nu, j})
\not \subset M$ and hence by Proposition~\ref{p}, it follows that
$\h f$ extends as an analytic set across the origin.

To complete the proof, it suffices to show that extension as an
analytic set implies extension as a locally proper holomorphic
correspondence. This is achieved in the next lemma.}
\end{case}

\begin{lem}
There exist neighbourhoods $U$ of $0$ and $U'$ of $0'$ such that
$F\subset U\times U'$ is a proper holomorphic correspondence which
extends $\hat f$.
\end{lem}

\begin{proof}
Extension as a holomorphic correspondence essentially follows from
\cite{[DP4]}. All nuances in the proof of Proposition~2.4 in
\cite{[DP4]} work in this situation as well provided the following
two modifications are made. Let $U, U'$ be neighbourhoods of $0,
0'$ respectively and suppose that $F \subset U \ti U'$ extends $\h
f$ as an analytic set in $U \ti U'$. Then it needs to be checked
that $F \cap (U^+ \ti U') \not= \emptyset$ and that there exists a
sequence $\{z_{\nu}\} \in M$ converging to $0$ such that $\h f$
extends as a correspondence across each $z_{\nu}$.

Suppose that $F \cap (U^+ \ti U') = \emptyset$. In this case, the
proof of Proposition~3.1 (or even Proposition~4.1 in \cite{[SV]})
shows that $(0, 0')$ is in the envelope of holomorphy of $\ov
{U^{-}} \ti U'$. The coefficients $a_{\mu \nu}(z)$ in \eqref{a}
can be regarded as holomorphic functions on $U^- \ti U'$ (i.e.,
independent of the $z'$ variables) and thus each $a_{\mu \nu}(z)$
extends holomorphically across $(0, 0')$. This extension must be
independent of the $z'$ variables by the uniqueness theorem and
hence $a_{\mu \nu}(z)$ extends holomorphically across the origin.
This shows that $\h f$ extends as a holomorphic correspondence
across the origin. To show the existence of the sequence
$\{z_{\nu}\}$ claimed above, let $\pi\hbox{:}~F \ra U$ be the
natural projection and define
\begin{equation*}
A = \{ (z, z') \in F\hbox{:}\ {\rm dim} \; ( \pi^{-1}(z) )_{(z, z')} \geq
1 \},
\end{equation*}
where $( \pi^{-1}(z) )_{(z, z')}$ denotes the germ of the fiber
over $z$ at $(z, z')$. Then $A$ is an analytic subset of $F$, and
since $F$ contains the graph of $\h f$ over $U^-$, it follows that
the dimension of $A$ is at most $n - 1$. Since Lipschitz maps do
not increase Hausdorff dimension, it follows that the Hausdorff
dimension of $\pi(A)$ is at most $2n - 2$. Pick $p \in M \sm
\pi(A)$. The fiber $F \cap \pi^{-1}(p)$ is discrete and this shows
that $\h f$ extends as a holomorphic correspondence across $p$.

Finally, we show that $U'$ can be chosen so small that the
projection $\pi'\hbox{:}~F\to U'$ is also proper. Indeed, for
$z'\in M'$, ${\pi}^{\prime -1}(z')$ is an analytic subset of $F$.
Since $\pi$ is proper, it follows by Remmert's theorem that $\h
F^{-1}(z') = \pi\circ {\pi}^{\prime -1} (z')$ is an analytic set.
The invariance property of Segre varieties yields $\hat F(Q_z\cap
U) \subset Q'_{z'}$ for any $z\in \h F^{-1}(z')$. Since $M$ is of
finite type, the set $\cup_{z\in \h F^{-1}(z')}Q_z$ has Hausdorff
dimension $n$, and therefore cannot be mapped by $\h F$ into
$Q'_{z'}$ which has dimension $n-1$. This shows that projection
$\pi'$ has discrete fibers on $M'$. It follows from the
Cartan--Remmert theorem that there exists a neighbourhood $U'$ of
$M'$ such that $\pi'$ has only discrete fibers, and therefore the
projection $\pi'$ from $F$ to $U'$ will be proper.

This completes the proof of Theorem~1.1.\hfill $\Box$\vspace{-.5pc}
\end{proof}

\section{Preservation of strata}

Fix $p \in M$ and let $p'_1, p'_2, \dots, p'_k \in \h f(p) \subset
M'$. Choose neighbourhoods $U, U'$ of $p, p'_1$ respectively and
let $\h f_1\hbox{:}~U^- \ra U'$ be a component of $\h f$ such that
$\h f_1(p) = p'_1$. Then $\h f_1$ extends as a holomorphic
correspondence $F \subset U \ti U'$ and to prove Theorem~1.2, it
suffices to focus on $\h f_1$, which will henceforth be denoted by
$\h f$. The following two general observations can be made in this
situation. First, the branching locus $\h \sigma$ of $\h F$ is an
analytic set in $U$ and the finite-type assumption on $M$ shows
that the real dimension of $\h \sigma \cap M$ is at most $2n - 3$.
The branching locus of $\h f$ denoted by $\si$, is contained in
$\h \si \cap U^-$. Second, the invariance property of Segre
varieties in \cite{[DP1]}, \cite{[V]} shows that $\h F$, the
extended correspondence, preserves the two components $U^{\pm}$.
That is, after possibly re-labelling ${U}^{\prime \pm}$, it
follows that $\h F(U^{\pm}) \subset {U}^{\prime \pm}$ and $\h F(M)
\subset M'$. The same holds for\break $\h G := \h
F^{-1}\hbox{:}~U' \ra U$.\vspace{-.2pc}

\setcounter{section}{1} \setcounter{theore}{1}
\begin{pot}{\rm
Let $p \in M^+$ and suppose that $\{\ze'_j\} \in M'$ is a sequence
converging to $p'_1$ with the property that the Levi form ${\cl
L}_{\rho}$ restricted to the complex tangent space to $M$ at
$\ze'_j$ has at least one negative eigenvalue. Fix $\ze'_{j_0} \in
U'$ for some large $j_0$. By shifting $\ze'_{j_0}$ slightly, we
may assume that $\ze'_{j_0} \notin \h \si' \cup \h F(M^0 \cap U)$,
where $\h \si'$ is the branching locus of $\h G$, and at the same
time retain the property of having at least one negative
eigenvalue. Let $g_1$ be a locally biholomorphic branch of $\h G$
near $\ze'_{j_0}$. Then $g_1(\ze'_{j_0})$ is clearly a
pseudoconvex point and this contradicts the invariance of the Levi
form. This shows that $\h f(M^+) \subset {M}^{\prime +}$. The same
arguments show that $\h f(M^-) \subset {M}^{\prime -}$.

Let $p \in M^+ \cap M^0$ and suppose that $p'_1 \in M^{\prime
+}_s$. The extending correspondence $\h F\hbox{:}~U \ra U'$
satisfies the invariance property, namely $\h F(Q_w) \subset
Q'_{w'}$ for all $(w, w') \in (U \ti U') \cap {\rm Graph}(\h F)$.
But near $p'_1$, the Segre map $\lambda$ is injective and this
shows that $\h F$, and hence $\h f$, is a single valued, proper
holomorphic mapping, say $f\hbox{:}~U \ra U'$ with $f(p) = p'_1$.
Two observations can be made at this stage: first, $f$ cannot be
locally biholomorphic near $p$ due to the invariance of the Levi
form. Second, if $V_f \subset U$ is the branching locus of $f$
defined by the vanishing of the Jacobian determinant of $f$, then
$V_f$ intersects both $U^{\pm}$. Indeed, suppose that $V_f \cap
U^- = \emptyset$. Choose a branch of $f^{-1}$ near some fixed
point $a' \in {U}^{\prime -}$ and analytically continue it along
all paths in ${U}^{\prime -}$ to get a well-defined mapping, say
$g\hbox{:}~{U}^{\prime -} \ra U^-$. The analytic set $F \subset U
\ti U'$ extends $g$ as a correspondence and hence \cite{[DP2]} $g$
is a well-defined holomorphic mapping in $U'$ and this must be the
single valued inverse of $f$. Thus $f$ is locally biholomorphic
near $p$ and this is a contradiction. The same argument works to
show that $V_f$ must intersect $U^+$ as well. Note that $V_f \cap
M$ has real dimension at most $2n - 3$. If $p \in \G_1$, choose
$U$ so small that $M^0 \cap U \subset \G_1$. Then there exists
\hbox{$q \in \G_1\!\sm\!(V_f \cap M)$} near $p$, where $f$ is
locally biholomorphic. Thus $q$ is mapped locally
biholomorphically to $f(q)$ which is a strongly pseudoconvex point
and this is a contradiction. If $p \in \G_3$, then again we shrink
$U$ so that $M^0 \cap U \subset \G_3$ and $(M \cap U)\!\sm\!\G_3
\subset M^+_s$. Then $f$ is locally biholomorphic near all points
in $(M \cap U)\!\sm\!\G_3$ and therefore $V_f \cap U^-$ must
cluster only along $\G_3$. Since the CR dimension of $\G_3 = n - 3
< (n - 1) - 1$ which is one less than the dimension of $V_f$, it
follows (Theorem~18.5 in \cite{[C]}) that $\ov {V_f \cap U^-}$ is
a closed, analytic set in $U$. Thus $V_f \cap U^-$ has two
analytic continuations, namely $V_f$ and $\ov {V_f \cap U^-}$ and
therefore they must be the same. This shows that $V_f$ cannot
intersect $U^+$ which is a contradiction. The same argument works
if $p \in \G_4$, the only difference being that $\ov {V}_f \subset
\ov {U^-}$ is analytic because of Shiffman's theorem. Thus if $p
\in M^+ \cap M^0$, then $p'_1 \in {M}^{\prime +} \cap {M}^{\prime
0}$.

To study this further, suppose that $p \in M^+ \cap \G_1$ and
$p'_1 \in M^{\prime +} \cap \G'_2$. Choose $U, U'$ small enough so
that $M^0 \cap U \subset \G_1$ and $M^{\prime 0} \cap U' \subset
\G'_2$. Pick \hbox{$q \in \G_1\!\sm\!(\h \si \cap M)$}. Then $\h
f$ splits near $q$ into finitely many well-defined holomorphic
mappings each of which extends across $q$. Moving $q$ slightly, if
needed, on \hbox{$\G_1\!\sm\!(\h \si \cap M)$}, each of these
holomorphic mappings are even locally biholomorphic near $q$.
Working with one of these mappings, say $f_1$, it follows that
$f_1(q) \notin {M}^{\prime +}_s$ due to the invariance of the Levi
form. This means that $f_1(q) \in \G'_2$. In the same way, all
points in $\G_1$ that are sufficiently near $q$ are mapped locally
biholomorphically by $f_1$ to $\G'_2$. This cannot happen as
$\G'_2$ has strictly smaller dimension than $\G_1$. The same
argument shows that $p'_1 \notin \G'_3 \cup \G'_4$. Hence $p'_1
\in M^{\prime +} \cap \G'_1$.

Suppose that $p \in M^+ \cap \G_2$ and $p'_1 \in M^{\prime +} \cap
\G'_1$. Considering $\h f^{-1}\hbox{:}~U' \ra U$, the arguments
used in the preceding lines show that this cannot happen. The case
when $p'_1 \in \G'_4$ can be dealt with similarly. Now suppose
that $p'_1 \in \G'_3$. As always, $U, U'$ will be small enough so
that $M^0 \cap U \subset \G_2$ and $M^{\prime 0} \cap U' \subset
\G'_3$. The arguments used above show that the cluster set of
points in $M^+_s \cap U$ is contained in $M^{\prime +}_s \cap U'$
and hence $\h f$ splits into finitely well-defined mappings each
of which is locally biholomorphic near points in $M^+_s \cap U$.
This shows that the branching locus $\si \subset U^-$ of $\h f$
clusters only along $\G_2$. Then $\h F(\si)$ is an analytic set of
dimension $n - 1$ in ${U}^{\prime -}$. There are two cases to
consider: first, if $\h F(\si)$ clusters only along $\G'_3$, then
arguing as above, $\ov {\h F(\si)} \subset \ov {{U}^{\prime -}}$
is a closed, analytic set in $U'$. The strong disk theorem shows
that $p'_1$ is in the envelope of holomorphy of ${U}^{\prime -}$
and this is a contradiction. Second, if there are points in $\ov
{\h F(\si)} \cap M^{\prime +}_s$, this means that $(\ov {\h F(\h
\si)} \cap M') \cap \G'_3$ has real dimension at most $2n - 4$.
Pick $q' \in \G'_3\!\sm\!({\ov {\h F(\h \si)}} \cap M')$ and note
that the continuity of $\h f$ implies that $\h f^{-1}(q') \in
M^+_s$. As seen above, this cannot happen. Thus $p'_1 \in \G'_2$
or $M^{\prime +}_s$. Similar arguments show that if $p \in M^+
\cap \G_3$ or $M^+ \cap \G_4$, then $p'_1 \in M^{\prime +} \cap
\G'_3$ or $M^{\prime +} \cap \G'_4$ respectively.

By reversing the roles of $U^{\pm}$, the same arguments used in
the preceding paragraphs can be applied to show that $\h f(M^-
\cap M^0) \subset {M}^{\prime -} \cap M^{\prime 0}$ with the
preservation of $M^- \cap \G_j$ for $j = 1, 3, 4$.

Finally, fix integers $i, j$ both at least $1$ such that $i + j =
n - 1$ and suppose that $p \in M_{i, j}$. Then there exists a
point $p_0$, in $U$ (chosen so small that $M \cap U \subset M_{i,
j}$) and arbitrarily close to $p$, where all branches of $\h f$
are well-defined and locally biholomorphic. The invariance of the
Levi form shows that the images of $p_0$ under any of the branches
of $\h f$ should all be in $M_{i, j}$. Note that each of these
images is close to $p'_1$. This cannot happen if $p'_1$ is in
${M}^{\prime +}, {M}^{\prime -}$ or in $M'_{i', j'}$ for $i \not=
i'$ and $j \not= j'$. The only possibility is that $p'_1$ is in
the relative interior of $\ov M'_{i, j}$. The same argument works
if $p$ is in the relative interior of $\ov M_{i, j}$.}\hfill
$\Box$
\end{pot}

\section*{Acknowledgements}

RS was supported in part by the Natural Sciences and Engineering
Research Council of Canada. KV was supported by DST (India) Grant
No.:~SR/S4/MS-283/05.

\end{document}